\documentclass[a4paper,11pt]{article}

\usepackage{jheppub} 
\usepackage{amsmath}
\usepackage[T1]{fontenc} 

\title{\boldmath Toy model of particle scattering theory}

\author[a]{Kyeong Min Kim}

\affiliation[a]{Mathematical Institute, University of Oxford, Woodstock Road, Oxford, OX2 6GG, UK}

\emailAdd{kyeong.kim@st-hughs.ox.ac.uk}

\abstract{The one dimensional probabilistic toy model of particle scattering theory is proposed. The toy model version of scattering probability is proved to be equal to the hypervolume of a $n$-dimensional figure. The solution for any $n$-particle toy model is presented as a contour integral, through Mellin trasnformation. The method of solving the contour integral is discussed. A nontrivial symmetry of this toy model, the invariance on initial position of the particles, is observed.}

\begin{document} 
\maketitle
\flushbottom

\section{Introduction}

Particle scattering is one of the greatest concern in physics. One of the main objectives in particle scattering is the calculation of the S-matrix, which is first brought by Wheeler \cite{a} and developed by Heisenberg \cite{b}. Some modern techniques include spinor helicity formalisms \cite{c} and bootstrapping the amplitude \cite{d}. However, calculating the scattering becomes extremely complicated when the number of the colliding particles gets large. This paper will suggest a simple toy model of particle scattering which allows calculation in a relatively elementary way.

We start building the toy model by confining the motion of particles in a single infinite line or in one dimension. Furthermore, we assume that the particles are of \textit{unit mass} for building the simpler model.

We label the particles as type A and type B according to the direction of the particle. Particles of type A move with constant speeds along the line in a positive direction, and particles of type B in a negative direction. We assume that the type B particles are initially all to the right of all the A particles. We can imagine this toy model as two particle accelerators in one dimension shooting beams of particle to each other.

Now, we imagine what would happen when two particles collide. Because the motion of the particles is confined in one dimension, deflecting with an angle would not be the case. In QFT, we encounter particle annihilation and creation \cite{e}, but we only consider annihilation, or possibly interpreted as absorption, of particles in this toy model. For simplicity, we propose a model with two kinds of outcomes after collision.

\begin{itemize}
    \item Different type of particle collides. Or two particles moving in opposite direction collide. One of the two particles is annihilated and the surviving particle continues to move with unchanged speed in the same direction as before.
    
    \item Same type of particle collides. Or two particles moving in same direction collide. The momentum of each particle does not change. They can move freely through each other.
    
\end{itemize}

The most fundamental nature of the particle scattering is violated in this model. The momentum is not conserved when different type of particles collides. However, we can take account of the momentum conservation by making the collision probabilistic. If the speeds of the type A and B particles are $a$ and $b$ respectively, then the probability that the type A particle survives is

\begin{equation}
    P(\text{A survives}) = \frac{a}{a+b}.
\end{equation}
 
 And the probability that the type B particle survive is, of course,
 
 \begin{equation}
 	P(\text{B survives}) = \frac{b}{a+b}.
 \end{equation}
 
 The probability that each particle survives is proportional to its velocity, which matches with our intuition that the particle with greater momentum should be more influential in the scattering process. Although the momentum is not conserved in individual collision, it is conserved on the average over many collisions. We can show that the expected value of the momentum after the collision is trivially conserved.
\begin{equation*}
<p_{\text{after}}> = \sum_{i} P(i \text{ survives}) p(i) = a-b.
\end{equation*}
The model that we built is a probabilistic model, which allows the momentum conservation, and directly follows our intuition that the probability to annihilate each other is proportional to its momentum.

We will consider the situation where a finite number $m$ of type A particles are in motion with speeds $(a_1, a_2,\dotsb,a_m)$ so as to encounter a finite number $n$ of type B particles moving in the opposite direction with speeds $(b_1,b_2,\dotsb,b_n)$. We will label each particles by $(A_1, A_2,\dotsb,A_m)$ and $(B_1,B_2,\dotsb,B_n)$ respectively. Annihilation process will continue until only one type of particles survive.

This paper focuses on the probability that the type B particle is entirely annihilated, calling this event \textit{A wins}. Likewise if the type A particles are all annihilated, we call this event \textit{B wins}. We will explore the methods of calculating the probability of \textit{A wins}.

In the following sections, we will manually calculate the probability of \textit{A wins} in simple cases, where $m$ and $n$ is a relatively small integer, and try to observe the property of the toy model by inspection. Next, we will try to look for the generalized formula of the probability of \textit{A wins}, as the hypervolume of the subset of the unit hypercube.

Solution as the hypervolume is still difficult to calculate. So we will transform the general solution in more solvable form, into the contour integral. Then we apply the solution using the contour integral into different cases of particle scattering. And we will look for the methods to solve the contour integrals.

\section{Elementary study of the toy model}
We will start with the simplest case where either $m$ or $n$ is 1. If $m=n=1$, then
\begin{equation*}
	\begin{split}
		&P(\text{A wins})=\frac{a_1}{a_1+b_1}.\\
		&P(\text{B wins})=\frac{b_1}{a_1+b_1}.
	\end{split}
\end{equation*}

This directly follows from the equation (1.1) and (1.2). If there are $m=1$ type A particle and $n$ type B particles, then the probability of the event \textit{A wins} is
\begin{equation*}
\begin{split}
	P(\text{A wins})&=P(A_1 \text{ annihilates } B_1) \, P(A_1 \text{ annihilates } B_1) \dotsb P(A_1 \text{ annihilates } B_1)\\
	& =\frac{{a_1}^2}{(a_1+b_1)(a_1+b_2)\dotsb(a_1+b_n)}.
\end{split}	
\end{equation*}
It is interesting that the order of $(b_1,b_2,\dotsb,b_n)$ does not affect the probability of A wins. It is true because it is just changing the order of multiplication in the equation. However, if we consider this as an actual particle scattering, this equation shows that the order in which A particle encounters the $n$ B particles does not matter. We shall carefully consider if the invariance of $P(\text{A wins})$ under changing the order of collision is just a coincidence in this simple case or a general property of this toy model. We shall look for the invariance in different cases as well.

Now let us consider another simple case where $m=n=2$. In this case, we can have two or three collisions before a type of particles is entirely annihilated. There can be many different orders in which collisions can occur, but for simplicity let us assume that particle $A_2$ and $B_2$ first collides and the surviving particle collides with left over particle. Then the probability of A wins is
\begin{equation*}
	\begin{split}
		P(\text{A wins})&=P(A_2\text{ wins }B_1) \big\{ 1-P(B_2 \text{ wins } A_1 \,\&\, A_2) \big\} \\
		&+ P(B_1\text{ wins }A_2) \big\{P(A_2 \text{ wins } B_1 \,\&\, B_2) \big\} \\
		&=\frac{{a_1}^2{a_2}^2+{a_1}^2a_2b_2+a_1{a_2}^2b_1+a_1a_2b_1b_2+{a_2}^2b_1b_2+{a_1}^2a_2b_1+{a_1}^2b_1b_2}{(a_1+b_1)(a_1+b_2)(a_2+b_1)(a_2+b_2)}.
	\end{split}
\end{equation*}
We can consider if $P(\text{A wins})$ is invariant under changing the order of collisions. The equation is symmetric about $a_1, a_2$ and $b_1, b_2$. So interchanging the particles does not affect the probability. The probability is still invariant, but we shall examine the general solution of this problem to see if the probability is invariant under particle exchange in this model.

Also, we can insert actual numbers in this problem, and observe what will happen. Let us examine two groups of particles, which have two particles in each group, with similar strength. Let $(a_1,a_2)=(30,20)$ and $(b_1,b_2)=(15,36)$. Inserting this into the equation, we obtain $P(\text{A wins})=270/539$, which is slightly greater than 0.5. The result is as anticipated.

Now, let us try obtaining the general solution by manual calculation. If there are $m$ type A particles and. $n$ type B particles, then the probability of A wins is (Let us assume that $A_1$ and $B_1$ first collides.)
\begin{equation*}
	\begin{split}
		&P(\text{A wins }m\text{ vs }n) \\
		&=P(B_1\text{ wins }A_1)\,P(A\text{ wins }m-1\text{ vs }n)+P(A_1\text{ wins }B_1)\,P(A\text{ wins }m\text{ vs }n-1)
	\end{split}
\end{equation*}
where
\begin{itemize}
	\item $P(\text{A wins }m-1\text{ vs }n) =P(B_1\text{ wins }A_{2})\,P(A\text{ wins }m-2\text{ vs }n) \\ +P(A_{2}\text{ wins }B_1)\,P(A\text{ wins }m-1\text{ vs }n-1)$
	\item $P(\text{A wins }m\text{ vs }n-1) =P(B_2\text{ wins }A_{1})\,P(A\text{ wins }m-1\text{ vs }n-1) \\ +P(A_{1}\text{ wins }B_2)\,P(A\text{ wins }m\text{ vs }n-2)$ 
	\\
	\vdots
\end{itemize}

It is still difficult to calculate the probability of A wins if $m$ and $n$ are large. But we do notice that the solution satisfies the recursive relation. It is actually the direct consequence of the conditional probability, where 
\begin{equation*}
	P(\text{A wins} \mid A_1 \text{ annihilates } B_1) = \frac{P(\text{A wins } \cap A_1 \text{ annihilates } B_1)}{P(A_1 \text{ annihilates } B_1)}
\end{equation*}
\begin{equation*}
	P(\text{A wins} \mid B_1 \text{ annihilates } A_1) = \frac{P(\text{A wins } \cap B_1 \text{ annihilates } A_1)}{P(B_1 \text{ annihilates } A_1)}
\end{equation*} 
Therefore, the recursive relation is
\begin{equation}
\begin{split}
	& P(\text{A wins }m\text{ vs }n) \\
	& = \frac{a_m}{a_m+b_1}P(\text{A wins} \mid A_1 \text{ annihilates } B_1)+\frac{b_1}{a_m+b_1}P(\text{A wins} \mid B_1 \text{ annihilates } A_1) \\
	& = \frac{a_m}{a_m+b_1}P(A\text{ wins }m-1\text{ vs }n)+\frac{b_1}{a_m+b_1}P(A\text{ wins }m\text{ vs }n-1)
\end{split}
\end{equation}
This recursive relation can be used to deduce the general solution.

\section{Solution as the hypervolume of the subset of the unit hypercube}
In this section, we will show that the probability of A wins is the hypervolume of the subset of the $m+n$ dimensional unit hypercube defined by condition:
\begin{equation}
	{x_1}^{a_1}{x_2}^{a_2}\dotsb {x_m}^{a_m} < {y_1}^{b_1}{y_2}^{b_2}\dotsb{y_n}^{b_n}
\end{equation}
or
\begin{equation}
	P(\text{A wins})=\int\dotsb \int_R dx_1 \,dx_2\dotsb dx_m \, dy_1 \,dy_2\dotsb dy_n
\end{equation}
\begin{equation*}
	\begin{split}
		\text{where }R & =\big\{ (x_1, \dotsb , x_m,y_1,\dots,y_n) \mid  0 \leq x_i \leq 1 \text{ and } 0\leq y_j \leq 1 \text{ for all } i \text{ and } j \big\} \\
& \cap \big\{ (x_1, \dotsb , x_m,y_1,\dots,y_n) \mid {x_1}^{a_1}{x_2}^{a_2}\dotsb {x_m}^{a_m} < {y_1}^{b_1}{y_2}^{b_2}\dotsb{y_n}^{b_n} \big\}.
	\end{split}
\end{equation*}

Before we move on to the proper proof, we should examine if this solution is invariant under particle exchange. If the order of $(a_1, a_2,\dotsb,a_m)$ and $(b_1,b_2,\dotsb,b_n)$ is changed, then we can find the corresponding subset of the unit hypercube just by reflecting or rotating the hypercube. Because the hypervolume is invariant about reflection and rotation, the probability is unchanged.

We should consider the following lemma in order to prove that the hypervolume gives the solution.

\textbf{Lemma 3.1} \textit{For any positive numbers $p$ and $q$,}
\begin{equation}
	\iint_R dx\,dy = \frac{a}{a+b} \int_S d\xi + \frac{b}{a+b} \int_T d\eta
\end{equation}
\begin{flushleft}
\textit{where $R$ is the subset of the unit square given by $px^a<qy^b$, $S$ is the subset of the unit interval given by $px^a < q$, and T is the subset of the unit interval given by $p<qy^b$.
}	
\end{flushleft}
\textit{Proof.} We change the variables. First, divide the integrating region by
\begin{equation*}
	R= \left( R \cap \{ (x,y) \mid x^a < y^b \} \right) \cup \left( R \cap \{ (x,y) \mid x^a \geq y^b \} \right)
\end{equation*}
where we define the first intersection as $R_1$ and the second intersection as $R_2$. Now we have
\begin{equation*}
	\iint_{R} dx\,dy=\iint_{R_1} dx\,dy+\iint_{R_2} dx\,dy
\end{equation*}
For the first term in the right hand side, let
\begin{equation*}
	x = \xi y^{b/a} \text{ or } \xi = \frac{x}{y^{b/a}}
\end{equation*}
We change the variable as $(x,y)$ to $(\xi,y).$ Then the corresponding Jacobian determinant is $J = y^{b/a}$. And therefore
\begin{equation*}
	\iint_{R_1} dx\,dy=\int y^{b/a} dy \int_{S} d\xi = \frac{a}{a+b} \int_S d\xi
\end{equation*}
where $S= \{ 0<\xi<1 \} \cap \{ p \xi^a < q \}$. Similarly, for the integrating region $R_2$, let 
\begin{equation*}
	y = \eta x^{a/b} \text{ or } \eta = \frac{y}{x^{a/b}}
\end{equation*}
We change the variable as $(x,y)$ to $(x, \eta)$. The Jacobian determinant is $J=x^{a/b}$. Therefore,
\begin{equation*}
	\iint_{R_2} dx \, dy = \int x^{a/b} dx \int_{T} d\eta =\frac{b}{a+b} \int_T d\eta
\end{equation*}
where $T=\{ 0<\eta<1 \} \cap \{ p < q \eta^b \}$. Therefore, 
\begin{equation*}
	\iint_R dx\,dy =\iint_{R_1} dx\,dy+\iint_{R_2} dx\,dy = \frac{a}{a+b} \int_S d\xi + \frac{b}{a+b} \int_T d\eta.
\end{equation*}
This is the end of the proof. Now let us apply the lemma with $x=x_1$, $y=y_1$, $a=a_1$, $b = b_1$. Then we have
\begin{equation*}
	\iint_R dx\,dy = \frac{a_1}{a_1+b_1} \int_S d\xi + \frac{b_1}{a_1+b_1} \int_T d\eta
\end{equation*}
where $R$ is the subset of the unit square given by ${x_1}^{a_1}\dotsb {x_m}^{a_m} < {y_1}^{b_1}\dotsb{y_n}^{b_n}$,
$S$ is the subset of the unit interval given by ${x_1}^{a_1}\dotsb {x_m}^{a_m} < {y_2}^{b_2}\dotsb{y_n}^{b_n}$,
and $T$ is the subset of the unit interval given by ${x_2}^{a_2}\dotsb {x_m}^{a_m} < {y_1}^{b_1}\dotsb{y_n}^{b_n}$. If we integrate the both sides of the upper equation about $x_2, \dotsb x_m, y_2, \dotsb , y_n$, we have
\begin{equation}
\begin{split}
	&\int \dotsi \int_{R'} dx_1 \dotsi dx_m \, dy_1 \dotsi dy_n \\
	&=\frac{a_1}{a_1+b_1} \int \dotsi \int_{S'} dx_1 \dotsi dx_m \, dy_1 \dotsi dy_n
	+\frac{b_1}{a_1+b_1} \int \dotsi \int_{T'} dx_1 \dotsi dx_m \, dy_1 \dotsi dy_n
\end{split}
\end{equation}
where $R'$ is the subset of the $(m+n)$ dimensional unit hypercube given by ${x_1}^{a_1}\dotsb {x_m}^{a_m} < {y_1}^{b_1}\dotsb{y_n}^{b_n}$,
$S'$ is the subset of the $(m+n-1)$ dimensional unit hypercube given by ${x_1}^{a_1}\dotsb {x_m}^{a_m} < {y_2}^{b_2}\dotsb{y_n}^{b_n}$,
and $T'$ is the subset of the $(m+n-1)$ dimensional unit hypercube given by ${x_2}^{a_2}\dotsb {x_m}^{a_m} < {y_1}^{b_1}\dotsb{y_n}^{b_n}$. We will define a function $f$ as (which is just the hypervolume of the subset)
\begin{equation}
	f(a_1,\dotsi,a_m;b_1,\dotsi, b_n)=\int \dotsi \int_{R'} dx_1 \dotsi dx_m \, dy_1 \dotsi dy_n.
\end{equation}
and the function $f$ satisfies the recurrence relation
\begin{equation}
\begin{split}
	&f(a_1,\dotsi,a_m;b_1,\dotsi, b_n)\\
	&=\frac{a_1}{a_1+b_1}f(a_1,\dotsi,a_m;b_2,\dotsi, b_n)
	+\frac{b_1}{a_1+b_1} f(a_2,\dotsi,a_m;b_1,\dotsi, b_n).
\end{split}
\end{equation}
The probability of A wins in the particle scattering model also satisfies the same recurrence relation as shown in the equation (2.1). Now we will prove that
\begin{equation}
	P(\text{A wins}) = f(a_1,\dotsi,a_m;b_1,\dotsi, b_n)
\end{equation}
using the proof by induction.

First, assume that the equation (3.7) is true for $m=m'-1, n=n'$ and $m=m', n=n'-1$. And assume that particle with speed $\underline{m}$ and $\underline{n}$ first collides. Then, 
\begin{equation*}
\begin{aligned}
		&P(\text{A wins};a_1,\dotsi,a_{m'};b_1,\dotsi,b_{n'})\\
		&= \frac{a_{\underline{m}}}{a_{\underline{m}}+b_{\underline{n}}}P(\text{A wins} \mid A_{\underline{m}} \text{ annihilates } B_{\underline{n}})+\frac{b_{\underline{n}}}{a_{\underline{m}}+b_{\underline{n}}}P(\text{A wins} \mid B_{\underline{n}} \text{ annihilates } A_{\underline{m}}) \\
		&=\frac{a_{\underline{m}}}{a_{\underline{m}}+b_{\underline{n}}}
		P(\text{A wins};a_1,\dotsi,a_{m'};b_1,\dotsi,b_{{\underline{n}}-1},b_{\underline{n}+1},\dotsi,b_{n'})\\
		&\hspace{2cm} +\frac{b_{\underline{n}}}{a_{\underline{m}}+b_{\underline{n}}}
		P(\text{A wins};a_1,\dotsi,a_{{\underline{m}}-1},a_{\underline{m}+1},\dotsi,a_{m'};b_1,\dotsi,b_{n'})\\
		&=\frac{a_{\underline{m}}}{a_{\underline{m}}+b_{\underline{n}}}
		f(a_1,\dotsi,a_{m'};b_1,\dotsi,b_{{\underline{n}}-1},b_{\underline{n}+1},\dotsi,b_{n'})\\
		&\hspace{2cm}+\frac{b_{\underline{n}}}{a_{\underline{m}}+b_{\underline{n}}}
		f(a_1,\dotsi,a_{{\underline{m}}-1},a_{\underline{m}+1},\dotsi,a_{m'};b_1,\dotsi,b_{n'})\\
		&=\frac{a_{\underline{m}}}{a_{\underline{m}}+b_{\underline{n}}}
		f(a_{\underline{m}},a_1,\dotsi,a_{{\underline{m}}-1},a_{\underline{m}+1},\dotsi,a_{m'};b_1,\dotsi,b_{{\underline{n}}-1},b_{\underline{n}+1},\dotsi,b_{n'})\\
		&\hspace{2cm}+\frac{b_{\underline{n}}}{a_{\underline{m}}+b_{\underline{n}}}
		f(a_1,\dotsi,a_{{\underline{m}}-1},a_{\underline{m}+1},\dotsi,a_{m'};b_{\underline{n}},b_1,\dotsi,b_{{\underline{n}}-1},b_{\underline{n}+1},\dotsi,b_{n'})\\
		&=f(a_{\underline{m}},a_1,\dotsi,a_{{\underline{m}}-1},a_{\underline{m}+1},\dotsi,a_{m'};b_{\underline{n}},b_1,\dotsi,b_{{\underline{n}}-1},b_{\underline{n}+1},\dotsi,b_{n'})\\
		&=f(a_1,\dotsi,a_{m'};b_1,\dotsi, b_{n'})
\end{aligned}
\end{equation*}
Therefore, the equation (3.0.7) holds for $m=m', n=n'$ if $m=m'-1, n=n'$ and $m=m', n=n'-1$ are true. Now, we should show that the statement holds for the initial values, where $m=1$ or $n=1$. For $m=1$ and $n=n'$,
\begin{equation*}
\begin{split}
		f(a_1;b_1,\dotsi,b_{n'})&=\int^{1}_{0}\dotsi\int^{1}_{0} \left( \int_{R'} dx_1 \right) dy_1\dotsi dy_{n'} \\
		&=\int^{1}_{0}{y_1}^{b_1/a_1}dy_1\dotsb \int^{1}_{0}{y_{n'}}^{b_{n'}/a_1}dy_{n'} \\
		&=\frac{{a_1}^{n'}}{(a_1+b_1)\dotsb(a_1+b_{n'})} \\
		&=P(\text{A wins};a_1;b_1,\dotsi,b_{n'})
\end{split}
\end{equation*}
Also, for $m=m'$ and $n=1$,
\begin{equation*}
	\begin{split}
		&f(a_1,\dotsi,a_{m'};b_1)\\
		&=\text{Hypervolume defined by ${x_1}^{a_1}\dotsb {x_{m'}}^{a_{m'}} < {y_1}^{b_1}$}\\
		&=1-\text{Hypervolume defined by ${y_1}^{b_1} < {x_1}^{a_1}\dotsb {x_{m'}}^{a_{m'}}$}\\
		&=1-f(b_1;a_1,\dotsi,a_{m'})\\
		&=1-\frac{{b_1}^{m'}}{(b_1+a_1)\dotsb(b_1+a_{m'})}\\
		&=P(\text{A wins}; a_1,\dotsi,a_{m'};b_1)
	\end{split}
\end{equation*}
By mathematical induction, the equation (3.7) is true for all natural number $m$ and $n$. And we have to extend the defintion of the fucntion $f$ such that
\begin{equation}
	f(a_1,\dotsi,a_m;)=1 \text{ and } f(;b_1,\dotsi,b_n)=0
\end{equation}
which is a coherent extension, because it means that the probability of A wins is 1 if there does not exist type B particle and the probability of A wins is 0 if there does not exist type A particle.

In this section, we have shown that the probability of A wins is equal to the hypervolume of the subset of unit hypercube. From this equality, we can easily observe the invariance of the probability under the exchange of the orders of particles, thus the initial position of the particles. Changing the order of the collision is equivalent to changing the axes of the hypercube or the order of multiplication of the terms in the function in the contour integral(which will be evident in the next section), which does not change the volume nor the integration.
\begin{align*}
\begin{split}
			P(\text{A wins})&=f(a_1,a_2,\dotsi,a_m;b_1,b_2,\dotsi,b_n) \\
			&=f(\underbrace{a_1',a_2',\dotsi,a_m'}_{\substack{\text{any permutation}\\\text{of } (a_1,\dotsi, a_m)}};\,\,\underbrace{b_1',b_2',\dotsi,b_n'}_{\substack{\text{any permutation}\\\text{of } (b_1,\dotsi, b_n)}}).
\end{split}
\end{align*}
However, solving the hypervolume is much more difficult than solving the probability manually. We will investigate for the methods of transforming the general solution into more solvable form.

\section{Transformation of the solution into the contour integral}
\subsection{Derivation of contour integral using the fourier transform}
We change the variable in the equation (3.2) by taking
\begin{equation*}
X_i=- \ln{x_i} \hspace{1cm} Y_j=-\ln{y_j}	
\end{equation*}
then we have
\begin{equation}
\begin{split}
&P(\text{A wins})\\
&=\int\dotsb \int_{R}\exp{\left( - \sum_i X_i - \sum_j Y_j \right)}dX_1 \,dX_2\dotsb dX_m \, dY_1 \,dY_2\dotsb dY_n
\end{split}
\end{equation}
where the region of integration $R$ is
\begin{equation*}
	\sum_{i=1}^{m} a_iX_i > \sum_{j=1}^{n} b_jY_j
\end{equation*}
for all $X_i, Y_j \geq 0$. Through the change of variable, the region of integration has become much simpler. Before we move on to the actual fourier transformation part, we should study the extension of the convolution theorem. Convolution is defined as
\begin{equation*}
	(f*g)(x) = \int_{-\infty}^{\infty}f(x-y)g(y)\,dy.
\end{equation*}
The convolution theorem states that the product of Fourier transform of functions $f$ and $g$ is equivalent to the Fourier transform of the convolution of $f$ and $g$ \cite{f}.
\begin{equation*}
\mathcal{F}(f*g)(t)= \widehat{f(t)}\widehat{g(t)}. 
\end{equation*}
And the extension of the convolution theorem is
\begin{equation}
	\mathcal{F}(f_1*f_2*\dotsi * f_n)(t)= \widehat{f_1(t)}\widehat{f_2(t)}\dotsb \widehat{f_n(t)}
\end{equation}
where
\begin{equation}
	(f_1*f_2*\dotsi * f_n)(x) = \int_{-\infty}^{\infty}\dotsi \int_{-\infty}^{\infty} f_1(x-x_2-\dotsi -x_n) f_2(x_2) \dotsi f_n(x_n) \,dx_2\dotsi dx_n.
\end{equation}
\textit{Proof of (4.1.2).} Substitute equation (4.1.3) to left hand side of (4.1.2).
\begin{equation*}
\begin{split}
	&\mathcal{F}(f_1*f_2*\dotsi * f_n)(t)\\
	&=\int_{-\infty}^{\infty} dx \,\exp{(-itx)}(f_1*f_2*\dotsi * f_n)(x) \\
	&=\int_{-\infty}^{\infty}\dotsi \int_{-\infty}^{\infty} f_1(x-x_2-\dotsi -x_n) f_2(x_2) \dotsi f_n(x_n) \,dx_2\dotsi dx_n e^{-itx}\,dx\\
	&=\int_{-\infty}^{\infty}\dotsi \int_{-\infty}^{\infty}\big\{ f_1(x-x_2-\dotsi -x_n) e^{-it(x-x_2-\dotsi-x_n)} \big\} \\ 
	& \hspace{4cm} \big\{ f_2(x_2) e^{-itx_2} \big\} \dotsi \big\{ f_2(x_2) e^{-itx_2} \big\}\,dx\,dx_2\dotsi dx_n
\end{split}
\end{equation*} 
Change variable such that $x-x_2-\dotsi-x_n=x_1$, we have
\begin{equation*}
\begin{split}
	&\mathcal{F}(f_1*f_2*\dotsi * f_n)(t)\\
	&=\int_{-\infty}^{\infty} f_1(x_1) e^{-itx_1} \, dx_1\int_{-\infty}^{\infty} f_2(x_2) e^{-itx_2} \, dx_2 \dotsm \int_{-\infty}^{\infty} f_n(x_n) e^{-itx_n} \, dx_n \\
	&=\widehat{f_1(t)}\widehat{f_2(t)} \dotsb \widehat{f_n(t)}.
	\end{split}
\end{equation*} 
Using this extension of the convolution theorem, we will transform the equation (4.1.1) to a contour integral on a complex plane. Change the variables such that $X_i=p_i/a_i$ and $Y_j=-q_j/b_j$. Then we can rewrite the equation (4.1.1) as
\begin{equation}
\begin{split}
		&P(\text{A wins})= {(a_1\dotsi a_m b_1 \dotsi b_n)}^{-1} \\
		& \int\dotsb \int_{R} \exp{\left( - \sum_i \frac{p_i}{a_i} + \sum_j \frac{q_j}{b_j} \right)}dp_1 \dotsb dp_m \, dq_1 \dotsb dq_n
\end{split}
\end{equation}
where $R$ is the region
\begin{equation*}
	\sum_{i=1}^{m} p_i + \sum_{j=1}^{n} q_j > 0, \hspace{0.5cm} p_i>0,\hspace{0.5cm} q_j<0.
\end{equation*}
We define piecewise functions $F_i$, $G_j$, and $H$ as
\begin{equation}
        F_i(x) =
        \begin{cases}
            \exp{(-\frac{x}{a_i})} &\hspace{1cm}x > 0 \\
            0 &\hspace{1cm} \text{otherwise}
        \end{cases}
\end{equation}
\begin{equation}
        G_j(x) =
        \begin{cases}
            \exp{(\frac{x}{b_j})} &\hspace{1cm}x < 0 \\
            0 &\hspace{1cm} \text{otherwise}
        \end{cases}
\end{equation}
\begin{equation}
        H(x) =
        \begin{cases}
            1 &\hspace{1cm} x < 0 \\
            0 &\hspace{1cm} \text{otherwise}
        \end{cases}
\end{equation}
Then equation (4.1.4) can be written in terms of functions $F_i$, $G_j$, and $H$.
\begin{multline}
P(\text{A wins})= {(a_1\dotsi a_m b_1 \dotsi b_n)}^{-1} \\
		 \int_{-\infty}^{\infty}\dotsb \int_{-\infty}^{\infty} F_1(p_1)\dotsi F_m(p_m) G_1(q_1)G_n(y_n) H(k-p_1-\dotsi-p_m-q_1-\dotsi-p_n)\\
		 dp_1 \dotsb dp_m \, dq_1 \dotsb dq_n	
\end{multline}
where $k=0$. We notice that equation (4.1.8) takes account of the integrating region $R$ using piecewise functions $F_i$, $G_j$, and $H$. We can show equation (4.1.8) as the convolution of functions.
\begin{equation}
	P(\text{A wins})={(a_1\dotsi a_m b_1 \dotsi b_n)}^{-1} \left( F_1 * \dotsi * F_m *G_1 * \dotsi * G_n * H \right)(0).
\end{equation}
Now, we can evaluate $P(\text{A wins})$ using the extension of the convolution theorem.
\begin{equation}
	\begin{split}
		&P(\text{A wins})\\
		&={(a_1\dotsi a_m b_1 \dotsi b_n)}^{-1} {\mathcal{F}}^{-1} \Big[ \mathcal{F} \left( F_1 * \dotsi * F_m *G_1 * \dotsi * G_n * H \right) \Big](0) \\
		&={(a_1\dotsi a_m b_1 \dotsi b_n)}^{-1}  {\mathcal{F}}^{-1} \Big[ \widehat{F_1} \dotsi \widehat{F_m} \,\, \widehat{G_1} \dotsi \widehat{G_n} \widehat{H}   \Big](0).
	\end{split}
\end{equation}
However, $\widehat{H}(t)$ is not integrable. So we should get round this by considering $H(x)$ as the limit of the function $H_{\epsilon}(x)$, where
\begin{equation}
        H_{\epsilon}(x) =
        \begin{cases}
            \exp{(\epsilon x)} &\hspace{1cm} x < 0 \\
            0 &\hspace{1cm} \text{otherwise}
        \end{cases}
\end{equation}
Then, the Fourier transform of $H_{\epsilon}(x)$ is 
\begin{equation}
		\widehat{H_{\epsilon}}(t) = \frac{1}{\epsilon -it}.
\end{equation}
Therefore,
\begin{equation}
\begin{split}
	 &\widehat{F_1} \dotsi \widehat{F_m} \,\, \widehat{G_1} \dotsi \widehat{G_n} \widehat{H_{\epsilon}}(t) \\
	 &= (a_1\dotsi a_m b_1 \dotsi b_n) \prod_{j=1}^{m} \frac{1}{1+ia_jt} \,\, \prod_{k=1}^{n} \frac{1}{1-ib_kt} \, \frac{1}{\epsilon - it}
\end{split}
\end{equation}
We take the inverse Fourier transform, for argument at $x=0$. We obtain \cite{g}
\begin{equation}
\begin{split}
		&P(\text{A wins})\\
		&={(a_1\dotsi a_m b_1 \dotsi b_n)}^{-1}  {\mathcal{F}}^{-1} \Big[ \widehat{F_1} \dotsi \widehat{F_m} \,\, \widehat{G_1} \dotsi \widehat{G_n} \widehat{H}   \Big](0)\\
		&=-\frac{1}{2\pi i} \int_{-\infty}^{\infty}dt \,\, \prod_{j=1}^{m} \frac{1}{1+ia_jt} \,\, \prod_{k=1}^{n} \frac{1}{1-ib_kt} \, \frac{1}{t+i\epsilon}.
	\end{split}	
\end{equation}
Now, we will move the integration path of the equation (4.14) to a closed contour in a complex plane. Let us consider the contour integration along the rectangular contour $C$ with its vertex at $(x,y)=(R,0), (R,k),(-R,k),$ and $(-R,0)$. If we take the limit such that $R \rightarrow \infty$ and $k \rightarrow \infty$, integration along the path $(R,0)\rightarrow(R,k)\rightarrow(-R,k)\rightarrow (-R,0)$ becomes 0. The only non-zero contribution is the linear path $(-R,0) \rightarrow (R,0)$, which is exactly $P(\text{A wins})$. Therefore,
\begin{equation}
		P(\text{A wins})=-\frac{1}{2\pi i} \oint_{C}dt \,\, \prod_{j=1}^{m} \frac{1}{1+ia_jt} \,\, \prod_{k=1}^{n} \frac{1}{1-ib_kt} \, \frac{1}{t+i\epsilon}.
\end{equation}
Consider the limit $\epsilon \rightarrow 0$. The Fourier transform of $H(t)$ was not integrable, so we have gotten round by redefining as $H(t)=\lim_{\epsilon \rightarrow 0} H_{\epsilon}(t)$. Substituting $\epsilon = 0$ is valid way of taking the limit, because  the contour $C$ only surrounds $a$-poles. Therefore,
\begin{equation}
\begin{split}
			P(\text{A wins})&=-\frac{1}{2\pi i} \oint_{C}dt \,\, \prod_{j=1}^{m} \frac{1}{1+ia_jt} \,\, \prod_{k=1}^{n} \frac{1}{1-ib_kt} \, \frac{1}{t}\\
			&=-\sum_{j=1}^{m} \text{Residue}({ia_j}^{-1},f_a(t))
\end{split}
\end{equation}
\cite{h} where
\begin{equation}
	f_a(t)=\prod_{j=1}^{m} \frac{1}{1+ia_jt} \,\, \prod_{k=1}^{n} \frac{1}{1-ib_kt} \, \frac{1}{t}.
\end{equation}
Also, if we change the variable such that $it=w$,
\begin{equation}
			P(\text{A wins})=-\frac{1}{2\pi i} \oint_{S}dw \,\, \prod_{j=1}^{m} \frac{1}{1-a_jw} \,\, \prod_{k=1}^{n} \frac{1}{1+b_kw} \, \frac{1}{w}
\end{equation}
We have written the probability of the event A wins as the contour integration, which is now much simpler than manual calculation. The strength of the contour integration comes from the residue theorem. We just have to calculate the residues of the $a$-poles and add them in order to obtain $P$(A wins). Also, we can notice the invariance of $P$(A wins) under the exchange of orders of the particles, because changing the order of the particles is equivalent to changing the order of multiplication, which does not affect the integration.
\subsection{Direct proof using partial fraction decomposition}
The equation (4.16) can be directly proved through the mathematical induction. We will show that the equation (4.16) satisfies the same recurrence relation with the particle scattering problem and the values for initial cases are identical. Before we move on to the proof of equation (4.16) by mathematical induction, we prove a lemma.

\textbf{Lemma 4.2.1} \textit{For any real numbers $c_i$,}
\begin{equation}
 \frac{1}{2\pi i} \oint dt  \,\,\frac{1}{w} \prod_{i=1}^{m} \frac{1}{1+c_iw}=0
 \end{equation}
\begin{flushleft}
\textit{where the contour integral includes all $c$-poles and the origin, or the entire plane.}
\end{flushleft}
\textit{Proof.} For simplicity, rewrite the equation (4.19) as
\begin{equation*}
	\oint \prod_{i=1}^{m} \frac{1}{z-d_i} dz=0.
\end{equation*}
for $d_i \in \mathbb{R}$ and R> $\lvert{d_i}\rvert$. Because the contour can be arbitrary as long as it surrounds every $c$-poles and the origin, we can take the contour as the circle with radius $2R$ centered at the origin. Then, for any $z$ on the circle with radius $2R$, 
\begin{equation*}
	\lvert{z-d_i}\rvert>R \text{ or } \Big\lvert{\frac{1}{z-d_i}}\Big\rvert<\frac{1}{R}.
\end{equation*}
And therefore,
\begin{equation*}
		\Bigg\lvert{\prod_{i=1}^{m}\frac{1}{z-d_i}}\Bigg\rvert=\prod_{i=1}^{m}\Bigg\lvert{\frac{1}{z-d_i}}\Bigg\rvert < \prod_{i=1}^{m} \frac{1}{R} = \frac{1}{R^m}.
\end{equation*}
And we substitute this result into the contour integration.
\begin{equation*}
	\Bigg\lvert{\oint dz \, \prod_{i=1}^{m}\frac{1}{z-d_i}}\Bigg\rvert \leq \oint dz \, \Bigg\lvert{\prod_{i=1}^{m}\frac{1}{z-d_i}}\Bigg\rvert<\oint \frac{dz}{R^m} =\frac{2\pi}{R^{m-1}}.
\end{equation*}
We can take $R$ as large as we want, so if we take the limit of $R \rightarrow \infty$,
\begin{equation*}
	\Bigg\lvert{\oint dz \, \prod_{i=1}^{m}\frac{1}{z-d_i}}\Bigg\rvert \rightarrow0
\end{equation*}
thus the equation (4.19) is true.

We could also prove the previous lemma by changing the variable as $w=1/z$. Now we move on to the direct proof of the equation (4.16) by mathematical induction. We will first show that it satisfies same recurrence relation that is also satisfied in the particle scattering problem. Partial fraction decomposition gives
\begin{equation*}
\frac{1}{1-a_i i t}\frac{1}{1+b_jit}=	\frac{a_i}{a_i+b_j}\frac{1}{1-a_i i t}+\frac{b_j}{a_i+b_j}  \frac{1}{1+b_jit}.
\end{equation*}
We define the function $g$ as
\begin{equation*}
	g(a_1,\dotsi,a_m;b_1,\dotsi,b_n) =-\frac{1}{2\pi i} \oint_{S}dt \,\, \prod_{i=1}^{m} \frac{1}{1+ia_it} \,\, \prod_{j=1}^{n} \frac{1}{1-ib_jt} \, \frac{1}{t}.
\end{equation*}
Using the partial fraction decomposition, we can show that $f$ satisfies the recurrence relation.
\begin{equation*}
	\begin{split}
		&g(a_1,\dotsi,a_m;b_1,\dotsi,b_n)\\
		&=-\frac{a_1}{a_1+b_1}\frac{1}{2\pi i} \oint_{S}dt \,\, \prod_{i=1}^{m} \frac{1}{1+ia_it} \,\, \prod_{j=2}^{n} \frac{1}{1-ib_jt} \, \frac{1}{t} \\
		&\hspace{3cm}-\frac{b_1}{a_1+b_1}\frac{1}{2\pi i} \oint_{S}dt \,\, \prod_{i=2}^{m} \frac{1}{1+ia_it} \,\, \prod_{j=1}^{n} \frac{1}{1-ib_jt} \, \frac{1}{t} \\
		&=\frac{a_1}{a_1+b_1}g(a_1,\dotsi,a_m;b_2,\dotsi,b_n)+\frac{b_1}{a_1+b_1}g(a_2,\dotsi,a_m;b_1,\dotsi,b_n).
	\end{split}
\end{equation*}
Now we will prove that
\begin{equation}
	P(\text{A wins};a_1,\dotsi,a_m;b_1,\dotsi,b_n) = g(a_1,\dotsi,a_m;b_1,\dotsi,b_n)
\end{equation}
using the mathematical induction. The procedure of the proof is identical to the proof in section 3. Functions $g$ and $f$ have same characteristics (they are both invariant on particle order exchange) and satisfies identical recurrence relation. We should only check if the initial values of $P(\text{A wins})$ matches the corresponding values of $g(a_1,\dotsi,a_m;b_1,\dotsi,b_n)$. For $m=1$ and $n=n'$,
\begin{equation*}
	g(a_1;b_1,\dotsi,b_{n'})=-\frac{1}{2\pi i} \oint_{S}dt \,\, \frac{1}{1+ia_1t} \,\, \prod_{j=1}^{n'} \frac{1}{1-ib_jt} \, \frac{1}{t}=\frac{{a_1}^{n'}}{(a_1+b_1)\dotsi(a_1+b_{n'})}
\end{equation*}
which is equivalent to $P(\text{A wins};a_1;b_1,\dotsi,b_{n'})$. Also, before we show the equivalence for $n=1$ case, we should prove a simple lemma.

\textbf{Lemma 4.2.2} \textit{For any positive real numbers $a_i$ and $b_j$,}
\begin{equation}
-\frac{1}{2\pi i} \oint_{S+T}dt \,\, \prod_{i=1}^{m} \frac{1}{1+ia_it} \,\, \prod_{j=1}^{n} \frac{1}{1-ib_jt} \, \frac{1}{t}
=1
 \end{equation}
\begin{flushleft}
\textit{where the contour $S$ surrounds the $a$-poles and $T$ surrounds the $b$-poles.}
\end{flushleft}
\textit{Proof.} From the Lemma 4.2.1, we directly have
\begin{equation*}
-\frac{1}{2\pi i} \oint_{S+T}dt \,\, \prod_{i=1}^{m} \frac{1}{1+ia_it} \,\, \prod_{j=1}^{n} \frac{1}{1-ib_jt} \, \frac{1}{t}=\frac{1}{2\pi i} \oint_{O}dt \,\, \prod_{i=1}^{m} \frac{1}{1+ia_it} \,\, \prod_{j=1}^{n} \frac{1}{1-ib_jt} \, \frac{1}{t}
 \end{equation*}
 where $O$ is the contour that only surrounds the origin. The contour integration around the origin is simply the residue at $t=0$. Therefore, the right hand side of the upper equation is $1$.
 Now, we should check the initial values for $m=m'$ and $n=1$.
 \begin{equation*}
 \begin{split}
	g(a_1,\dotsi,a_{m'};b_1)&=-\frac{1}{2\pi i} \oint_{S}dt \,\, \prod_{i=1}^{m'} \frac{1}{1+ia_it} \,\, \frac{1}{1-ib_1t} \, \frac{1}{t}\\
	&= 1+\frac{1}{2\pi i} \oint_{T}dt \,\, \prod_{i=1}^{m'} \frac{1}{1+ia_it} \,\, \frac{1}{1-ib_1t} \, \frac{1}{t}\\
	&= 1-\frac{{b_1}^{m'}}{(b_1+a_1)\dotsi(b_1+a_{m'})}\\
	&=P(\text{A wins};a_1,\dotsi,a_{m'};b_1).
 \end{split}
 \end{equation*}
By mathematical induction, the equation (4.20) holds for every natural number $m$ and $n$. The significance of the integration about the contour $T$ is that it is $P(\text{B wins})$. We can directly show this equality through changing the variable as $t=-t'$. Then everything is identical with the equation (4.16) except that the positions of $a_i$ and $b_j$ are changed. So changing the variable as $t=-t'$ is identical with the reflection of the particle scattering system, such that $B$ particles are moving in positive direction and $A$ particles are moving in negative direction. Also, Lemma 4.2.2 suggests that
\begin{equation}
	P(\text{A wins})+P(\text{B wins})=1
\end{equation}
which seems to be obvious, because the particle scattering only results in two independent events A wins and B wins.
\section{Application of the contour integral}
\subsection{Distinct case}
If all $a_i$ are distinct, the general solution becomes very simple.
\begin{equation}
\begin{split}
				P(\text{A wins})=-\sum_{i=1}^{m} \text{Residue}({a_i}^{-1},f_a(w)).
\end{split}
\end{equation}
where
\begin{equation}
	f_a(w)=\prod_{i=1}^{m} \frac{1}{1-a_iw} \,\, \prod_{j=1}^{n} \frac{1}{1+b_jw} \, \frac{1}{w}.
\end{equation}
If $a_i \neq a_k$ for $i \neq k$,
\begin{equation}
\begin{split}
				P(\text{A wins})&=\sum_{i=1}^{m}\prod_{\substack{k=1\\ k\neq i}}^{m} \frac{a_i}{a_i-a_k} \,\, \prod_{j=1}^{n} \frac{a_i}{a_i+b_j}
				 \\
				&=\sum_{i=1}^{m+n}\prod_{\substack{k=1\\ k\neq i}}^{m} \frac{a_i}{a_i-a_k} .
\end{split}
\end{equation}
where $a_{m+r}=-b_r$.
\subsection{Non-distinct case}
If some $a_i$ are not distinct, we cannot use the simplified formula such as the equation (5.3). We have calculated the residues by assuming that the $a$-poles are simple poles. The calculations become complicated when we have to calculate the residues of poles with order higher than 1. 
\subsubsection{All particles with same speed}
Consider a simple case where a set of $m$ particles, all of speed 1, collide with a set of $n$ particles, also with speed 1. Then the probability of A wins is
\begin{equation}
	P(\text{A wins}) = -\frac{1}{2\pi i} \oint_S \frac{1}{t} {\left( \frac{1}{1-t} \right)}^m  {\left( \frac{1}{1+t} \right)}^n \, dt = s(m,n). 
\end{equation}
Before we move on to actual calculation, we can infer the properties of $s(m,n)$, considering that the speed of every particles in the system is identical.
\begin{equation*}
	\begin{cases}
		s(m,0)=1\\
		s(0,n)=0\\
		s(m,n)=1-s(n,m)\\
		s(k,k)=\frac{1}{2}
	\end{cases}
\end{equation*}
First three properties are general properties of $P(\text{A wins})$. However, the last property is unique. We should check this property after calculating the probability.
\begin{equation}
	\begin{split}
		P(\text{A wins}) & = \frac{{-1}^{m+1}}{2\pi i} \oint_S \frac{1}{t} {\left( \frac{1}{t-1} \right)}^m  {\left( \frac{1}{1+t} \right)}^n \, dt\\
		&=\frac{{-1}^{m+1}}{(m-1)!}\frac{d^{m-1}}{{dt}^{m-1}} \left( \frac{1}{t{(1+t)}^n} \right)\Biggr\rvert_{t = 1}
		=\sum_{i=0}^{m-1}{{n+i-1}\choose{i}}2^{-n-i}.
	\end{split}
\end{equation}
However, the last property of the solution $s(k,k)=\frac{1}{2}$ does not seem to be obvious in the equation (5.5). We can directly calculate the equation (5.4) to see if $P(\text{A wins})=0.5$ if $m=n$. Consider a rectangular contour $C_r$ with vertex at $(R,0), (R,ih), (-R,ih),$ and $(-R,0)$ with a semicircle with radius $r$ in the origin. We take the limit of $R \rightarrow \infty$, $h \rightarrow \infty$, and $r\rightarrow0$. In this limit, the integration along the path $(R,0)\rightarrow (R,ih)\rightarrow (-R,ih)\rightarrow (-R,0)$ becomes $0$. Also, adding the integration along the path $(-R,0) \rightarrow (-r,0)$ and $(r,0) \rightarrow (R,0)$ gives 0, because the integrand is an odd function. The only path that gives non zero result is the semi circular path. Therefore,
\begin{equation*}
	s(n,n)= -\frac{1}{2\pi i} \oint_S \frac{1}{t} {\left( \frac{1}{1+t^2} \right)}^n \, dt =\lim_{r\rightarrow0} -\frac{1}{2\pi i}\int_{\pi}^{0}\frac{d\theta}{re^{i\theta}(1+r^2e^{2i\theta})^n}=\frac{1}{2}.
\end{equation*}
\subsubsection{Each type particles with same speed}

Let us think about another simple case, where a set of $m$ particles, all of speed 1, collide with a set of $n$ particles, all of speed $v$. Then the probability of A wins is
\begin{equation}
		P(\text{A wins}) = -\frac{1}{2\pi i} \oint_S \frac{1}{t} {\left( \frac{1}{1-t} \right)}^m  {\left( \frac{1}{1+vt} \right)}^n \, dt = s_v(m,n).
\end{equation}
The calculation is very analogous to the section 5.2.1, which results
\begin{equation}
	s_v(m,n)= \sum_{i=0}^{m-1} {{n+i-1}\choose i}\frac{v^i}{{(1+v)}^{n+i}}.
\end{equation}

\subsection{Matching and beating of the groups of the particles}
We define that two sets of particles are matched if both of them have probabilities of 0.5 of winning when colliding each other. We also define a set of particle beats another set of particle on average if the probability of winning is greater than 0.5. Matching has a counterintuitive property: matching is not a transitive relation. For example, we could consider a single particle of speed 60 and two-particle set $(20,30)$, $(15,36)$, and $(12,40)$.
\begin{equation*}
	\begin{cases}
		P(\text{A wins};60;20,30)= \frac{1}{2} \\
		P(\text{A wins};60;15,36)= \frac{1}{2} \\
		P(\text{A wins};60;12,40)= \frac{1}{2} \\
		P(\text{A wins};20,30;15,36)= \frac{270}{539} >\frac{1}{2} \\
		P(\text{A wins};15,36;12,40)= \frac{314}{627} >\frac{1}{2}\\
		P(\text{A wins};12,40;20,30)= \frac{293}{588}<\frac{1}{2}
	\end{cases}
\end{equation*}
From the above example, it seems that beating on average is a transitive relation. $(20,30)$ beats $(15,36)$, $(15,36)$ beats $(12,40)$, and $(20,30)$ beats $(12,40)$. However, beating on average is not a transitive relation. There exist sets of particles $P,Q,R$ such that on average $P$ beats $Q$, $Q$ beats $R$, but $R$ beats $Q$. We shall find the counterexample. For convenience, we set $Q$ as a single particle with speed 1, and $P, R$ as two-particle sets such that $P=(c,d)$ and $R=(a,b)$. We search for the solutions graphically, as in the figure 1, such that the group of particle $(x,y)$ above the curve beats $1$, and $1$ beats $(x,y)$ below the curve. And $(x,y)$ on the curve matches with $1$. As we have seen in the previous example, the group of two particles is stronger as the speeds are closer. $(20,30)$ is stronger than $(12,40)$, even though the net momentum is smaller. Therefore, we can take $R$ as the point near the middle of the curve, but right below and take $P$ as the point near the end of the curve, but right above. For example, beating on average is not a transitive relation for $(1)$, $(0.414213,0.414212)$, and $(0.9,0.0526317)$.
\begin{figure}[tbp]
\centering 
\includegraphics[width=.4\textwidth]{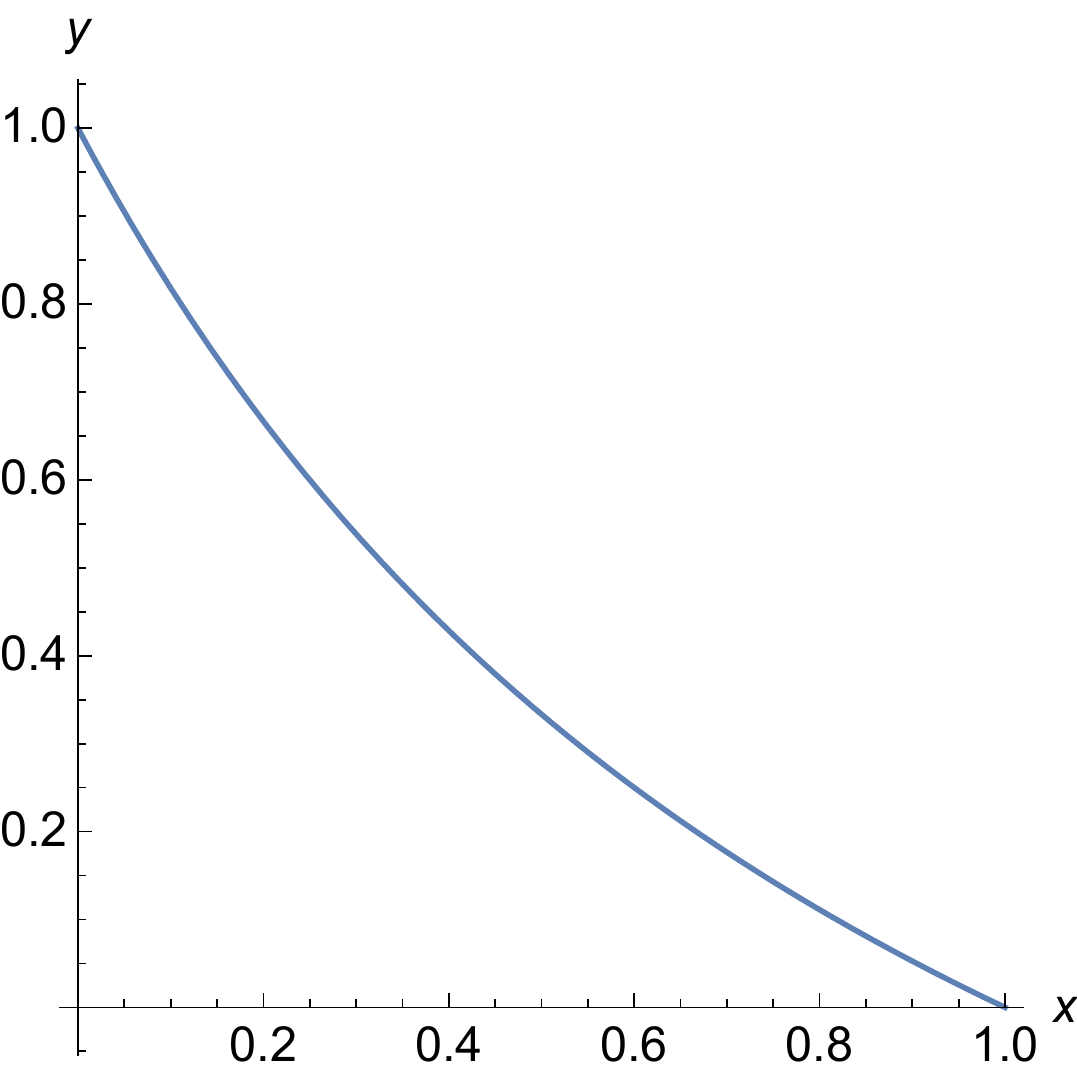}
\caption{\label{fig:1} $\frac{1}{1+x}\frac{1}{1+y}=\frac{1}{2}$.}
\end{figure}
\section{Solving the contour integral}
\subsection{Approximation method}
For the general particle scattering problem, the speed of the particles may not be distinct, and it may not be the simple cases where the speed of particles are all identical, as we have dealt in section 5.2.

We can approximate the probability of A wins by giving a slight variation to the non distinct speed. We can generalize the non distinct case by $x_1$ $A_1$-particles, $x_2$ $A_2$-particles, $\dotsb$, $x_m$ $A_m$-particles, $y_1$ $B_1$-particles, $y_2$ $B_2$-particles, $\dotsb$, and $y_n$ $B_n$-particles (where $a_1,\dotsi,a_m,b_1,\dotsi,b_n$ are distinct) such that
\begin{equation}
\begin{split}
		P(\text{A wins})&=f(\underbrace{a_1,\dotsi,a_1}_{x_1},\dotsi,\underbrace{a_m,\dotsi,a_m}_{x_m},\underbrace{b_1,\dotsi,b_1}_{y_1},\dotsi,\underbrace{b_n,\dotsi,b_n}_{y_n})\\
		&=-\frac{1}{2\pi i} \oint_{S}dt \,\, \prod_{i=1}^{m} \frac{1}{{(1+ia_it)}^{x_i}} \,\, \prod_{j=1}^{n} \frac{1}{{(1-ib_jt)}^{y_j}} \, \frac{1}{t}.
\end{split}
\end{equation}
We can select an arbitrary positive real number $\epsilon$ such that
\begin{equation*}
	\Delta a \gg \epsilon \text{ and } \Delta b \gg \epsilon 
\end{equation*}
where $\Delta a$ and $\Delta b$ is difference between any $a_i$ and $b_j$. We can give a slight variation to $x_i$ $a_i$-particles and $y_j$ $b_j$-particles as
\begin{equation*}
	\begin{split}
		&\text{1st }a_i \rightarrow a_i + \epsilon = a_{i,1}\\
		&\text{2nd }a_i \rightarrow a_i + 2\epsilon = a_{i,2}\\
		&\hspace{0.7cm} \vdots \\
		&x_i\text{th }a_i \rightarrow a_i + x_i\epsilon = a_{i,x_i}
	\end{split}
\end{equation*}
and
\begin{equation*}
	\begin{split}
		&\text{1st }b_j \rightarrow b_j + \epsilon = b_{j,1}\\
		&\text{2nd }b_j \rightarrow b_j + 2\epsilon = b_{j,2}\\
		&\hspace{0.7cm} \vdots \\
		&y_j\text{th }b_j \rightarrow b_j + x_i\epsilon = b_{j,y_j}.
	\end{split}
\end{equation*}
Now, all speeds are distinct and use the equation (5.1.4).
\begin{equation}
\begin{split}
		P(\text{A wins})\approx \sum_{i=1}^{m} \sum_{q=1}^{x_i} \prod_{\substack{u=1\\ u\neq q}}^{x_i} \frac{a_{i,q}}{a_{i,q}-a_{i,u}}\prod_{\substack{k=1\\k\neq i}}^{m}\prod_{s=1}^{x_k} \frac{a_{i,q}}{a_{i,q}-a_{k,s}}\prod_{j=1}^{n} \prod_{r=1}^{y_j} \frac{a_{i,q}}{a_{i,q}+b_{j,r}}\\
		=\sum_{i=1}^{m} \sum_{q=1}^{x_i} \prod_{\substack{u=1\\ u\neq q}}^{x_i} \frac{a_{i}+q\epsilon}{(q-u)\epsilon}\prod_{\substack{k=1\\k\neq i}}^{m}\prod_{s=1}^{x_k} \frac{a_{i}+q\epsilon}{a_{i}-a_{k}+(q-s)\epsilon}\prod_{j=1}^{n} \prod_{r=1}^{y_j} \frac{a_{i}+q\epsilon}{a_{i}+b_{j}+(q+r)\epsilon}.
\end{split}
\end{equation}
The approximation will be more accurate as we take smaller value of $\epsilon$.
\subsection{General Solution}
Previous section calculated the proability with the approximation method, but not the exact solution. In this section, we will calculate the exact solution using two different methods.
\subsubsection{Taking the limit of the approximation method}
We can obtain the exact solution when we take the limit of $\epsilon \rightarrow 0$ of the equation (6.2). After a lengthy calculation, we obtain
\begin{align}
	\begin{split}
		P(\text{A wins})&=\sum^{m}_{i=1}\frac{1}{(x_i-1)!\left( \prod^{m+n}_{\substack{k=1\\k\neq i}} (a_i-a_k)^{|x_k|} \right)} \\
		&\cdot \sum^{x_i}_{q=1}\left( \sum^{x_i-1}_{w=0} {{c}\choose{w}}{a_i}^{c-w}q^w(-1)^{x_i-q}{{x_i-1}\choose{q-1}} \left( \sum_{\mathcal{H}^{x_i-1-w}_{i,\,q}} \prod_{H^{x_i-1-w}_{i,\,q}(d)} \frac{q''-s}{a_i-a_k} \right)  \right)
	\end{split}
\end{align}
where
\begin{equation*}
	\sum^{m}_{k=1}x_k+\sum^{j=1}_{n}y_j-1=c
\end{equation*}
\begin{equation*}
	a_{m+j}=-b_j \text{ and } x_{m+j}=-y_j \text{ for } 1 \leq j \leq n.
\end{equation*}
$\mathcal{H}_{i,q}$ is defined as the set of all possible $(k,q'',s)$ such that
\begin{align*}
	\begin{split}
		\mathcal{H}_{i,q}=& \bigg\{ (k,q'',s) \mid \{q''\mid 1\leq q'' \leq x_i,q''\neq q, q'' \neq s\}, \{k \mid 1 \leq k \leq m+n , k \neq i\}, \\ 
		&\{ s \mid 1 \leq s \leq x_k \, \cup \, x_k \leq s \leq -1  \} \bigg\}.
		 \end{split}
\end{align*}
And $\mathcal{H}^z_{i,q}$ is defined as the set of all possible $z$-combination of the set $\mathcal{H}_{i,q}$. $H^{z}_{i,q}(d)$ is defined as the $d$th $z$-combination, such that $ H^{z}_{i,q}(d) \in \mathcal{H}^z_{i,q}$, and $1 \leq d \leq {{|\mathcal{H}_{i,q}|}\choose{z}}$, where $|\mathcal{H}_{i,q}|$ is the cardinality of the set $\mathcal{H}_{i,q}$. It can be easily shown that equation (6.3) simplifies to (5.3) for the distinct cases.

\subsubsection{Differentiation respect to parameters}

We can also differentiate respect to the parameters, $a_i$s, to obtain the general solution. We rewrite the equation (6.1) as
\begin{equation}
\begin{split}
		P(\text{A wins})&=-\frac{1}{2\pi i} \oint_{S}dt \,\, \prod_{i=1}^{m} \frac{1}{{(1-a_it)}^{x_i}} \,\, \prod_{j=1}^{n} \frac{1}{{(1+b_jt)}^{y_j}} \, \frac{1}{t} \\
		&=-\frac{1}{2\pi i} \oint_{S}dt \,\, \prod_{i=1}^{m} \left( \frac{1}{(x_i-1)!\,t^{x_i-1}}\frac{\partial^{x_i-1}}{\partial a^{x_i-1}} \left(\frac{1}{{1-a_it}}\right) \right) \,\, \prod_{j=1}^{n} \frac{1}{{(1+b_jt)}^{y_j}} \, \frac{1}{t} \\
		&= -\frac{1}{2\pi i}\prod_{l=1}^{m} \frac{1}{(x_l-1)!}\frac{\partial^{x_l-1}}{\partial a^{x_l-1}}
		 \oint_{S}dt \,\, \underbrace{t^{m-1-\sum x_i} \prod_{i=1}^{m} \frac{1}{1-a_it} \,\, \prod_{j=1}^{n} \frac{1}{{(1+b_jt)}^{y_j}}}_{=\,\Phi (t)} \\
		&= -\prod_{l=1}^{m} \frac{1}{(x_l-1)!}\frac{\partial^{x_l-1}}{\partial a^{x_l-1}}\left( \sum_{i=1}^{m} \text{Residues} \left( \Phi (t),\frac{1}{a_i} \right) \right) \\
		&=\prod_{l=1}^{m} \frac{1}{(x_l-1)!}\frac{\partial^{x_l-1}}{\partial a^{x_l-1}} \left( \sum_{i=1}^{m} a_i^{\sum x_i -m} \prod^{m}_{\substack{k=1\\k\neq i}} \frac{{a_i}}{a_i-a_k} \prod_{j=1}^{n} \frac{a_i}{a_i+b_j}
\right).
\end{split}
\end{equation}

\section{Conclusion}

The one dimensional toy model of particle scattering is constructed. Two different types of particles collide and annihilate the other particle with a probability proportional to their momentum. The surviving particle moves with same momentum. We name each type of particles as A and B, and we presume that same type of particles can pass through each other, without changes in momentum. The toy model reflects the physical nature through the momentum conservation on average. The collisions will continue until only a certain type of the particles remains, and we aim to calculate the probability of \textit{A wins} when the type B particle is entirely annihilated.

The hypervolume of the subset of the unit hypercube $R$ is proven to be equal to the probability of A wins using the mathematical induction. The hypervolume is expressed as the integral of products of piecewise functions, and it is transformed to a contour integral by using the convolution theorem and the Fourier transform. We then calculate the contour integral in different cases, such as when the speed of the type A particles are distinct, all particles have same speed, and each type of particles have same speed. The relation between groups of particle, matching and beating, is defined. Matching and beating are proven to be intransitive relations. 

Methods to solve the probability of A wins in the general case is studied. Approximation method is used to approximate the non distinct general case as distinct case, and the exact solution is obtained through taking the limit. Also, the contour integral is solved by differentiating respect to the parameters.

The most nontrivial property of this toy model is that the solution is invariant about the initial position of the particles, which is the direct result of the probability of A wins being equal to the hypervolume and the contour integral. The toy model failed to account for interactions that resemble the standard model, but we could find out how the contour integrals naturally arise from this physically inspired toy model, and explore basic mathematical techniques which can be applied in calculating probabilities.

\acknowledgments

I thank professor Andrew Hodges for supervising the project.


\end{document}